\documentclass[draft]{amsart}

\usepackage{amsthm}
\usepackage[leqno]{amsmath}
\usepackage{latexsym,amsfonts,amssymb}
\usepackage[all]{xy} \SelectTips{eu}{}
\usepackage{hyperref}
\usepackage{enumerate}
\usepackage{amssymb,amstext,amsmath,amsthm,amsfonts,enumerate,latexsym}

\newcommand{\numberseries}{\bfseries}   %Fontseries usced for numbering
\newlength{\thmtopspace}                %Space above theorem
\newlength{\thmbotspace}                %Space below theorem
\newlength{\thmheadspace}               %Space after theorem label
\newlength{\thmindent}                  %For indenting
\newtheoremstyle{bfupright head,upright body}
                {\thmtopspace}{\thmbotspace}{\upshape}
                {\thmindent}{\bfseries}{.}{\thmheadspace}
                {\thmnote{#3}{\numberseries \thmnumber{\:#2}}}

\newtheoremstyle{fixed bf head,slanted body}
                {\thmtopspace}{\thmbotspace}{\slshape}
                {\thmindent}{\bfseries}{.}{\thmheadspace}
                {\thmname{#1}{\numberseries \thmnumber{\:#2}}\thmnote{ (#3)}}

\newtheoremstyle{fixed bf head,upright body}
                {\thmtopspace}{\thmbotspace}{\upshape}
                {\thmindent}{\bfseries}{.}{\thmheadspace}
                {\thmname{#1}{\numberseries \thmnumber{\:#2}}\thmnote{ (#3)}}

\newtheoremstyle{unnumbered paragraph}
                {\thmtopspace}{\thmbotspace}{\upshape}
                {\parindent}{\upshape}{}{0pt}

\theoremstyle{fixed bf head,slanted body}
\newtheorem*{thm*}{Theorem}
\newtheorem{thm}{Theorem}
\newtheorem*{prp*}{Proposition}
\newtheorem{prp}[thm]{Proposition}
\newtheorem*{cor*}{Corollary}
\newtheorem{cor}[thm]{Corollary}
\newtheorem*{lem*}{Lemma}
\newtheorem{lem}[thm]{Lemma}
\theoremstyle{fixed bf head,upright body}
\newtheorem{dfn}[thm]{Definition}
\newtheorem{rmk}[thm]{Remark}

\theoremstyle{unnumbered paragraph}
\newtheorem{ipg*}{}

\theoremstyle{bfupright head,upright body}
\newtheorem*{bfhpg*}{}
\newtheorem{bfhpg}[thm]{}

\newenvironment{prf*}[1][Proof]{%
  \begin{proof}[\bf #1]
    \setcounter{equation}{0}
    }
  {\end{proof}
}

\renewcommand{\eqref}[1]{(\ref{eq:#1})}
\newcommand{\nt}[2][$^\diamondsuit$]{%
        \hspace{0pt}#1\marginpar{\tt\raggedleft #1 #2}}
\def\urltilda{\kern -.15em\lower .7ex\hbox{\~{}}\kern .04em} 
\newcommand{\NN}{\mathbb{N}}

\renewcommand{\nt}[2][]{} %Hides comments in tt   
\newcommand{\set}[2][\mspace{1mu}]{\{#1 #2 #1\}}
\newcommand{\setof}[3][\mspace{1mu}]{\{#1#2 \mid #3#1\}}
\newcommand{\Spec}[1]{\operatorname{Spec}#1}
\newcommand{\p}{\mathfrak{p}}

\newcommand{\q}{\mathfrak{q}}
\renewcommand{\H}[2][]{\operatorname{H}_{#1}(#2)}
\newcommand{\SingR}{\operatorname{Sing}R}
\newcommand{\hgt}[2][R]{\operatorname{ht}_{#1}#2}
\newcommand{\codim}[1]{\operatorname{codim}#1}
\newcommand{\lra}{\longrightarrow}
\newcommand{\Ker}[1]{\nobreak{\operatorname{Ker}#1}}
\renewcommand{\Im}[1]{\nobreak{\operatorname{Im}#1}}

\newcommand{\pd}[2][R]{\operatorname{pd}_{#1}#2}
\newcommand{\sing}{\mathsf{S}(R)}
\newcommand{\RegR}[1][R]{\operatorname{Reg}#1}
\renewcommand{\mod}[1]{\mathsf{mod}(#1)}

\newcommand{\tier}[2]{\mathsf{tier}_{#1}#2}
\renewcommand{\cos}[1]{\mathsf{cosyz}\,#1}
\newcommand{\close}[1]{\langle\mspace{1mu} #1 \mspace{1mu}\rangle}

\def\p{\mathfrak{p}}

\def\mod{\operatorname{\mathsf{mod}}}
\def\dim{\operatorname{dim}}
\def\add{\operatorname{\mathsf{add}}}
\def\depth{\operatorname{depth}}
\def\fl{\operatorname{\mathsf{fln}}}
\def\cosyz{\operatorname{\mathsf{cosyz}}}

\def\PD{\mathsf{fpd}}

\def\dep{\operatorname{\mathsf{dep}}}

\def\S{\mathsf{S}}
\def\T{\mathsf{T}}

\def\cm{\mathsf{CM}}
\def\Ext{\operatorname{Ext}}
\def\s{\mathsf{S}}

\begin{document}

\title{Building modules from the singular locus}

\author{Jesse Burke}

\address{%Department of Mathematics,
  University of California, Los Angeles, CA 90095, U.S.A.}

\email{jburke@math.ucla.edu}

\urladdr{http://www.math.ucla.edu/\urltilda jburke}

\author{Lars Winther Christensen}

\address{Texas Tech University, Lubbock, TX 79409, U.S.A.}

\email{lars.w.christensen@ttu.edu}

\urladdr{http://www.math.ttu.edu/\urltilda lchriste}

\author{Ryo Takahashi}

\address{%Graduate School of Mathematics,
  Nagoya University, Furocho, Chikusaku, Nagoya 464-8602, Japan}

\email{takahashi@math.nagoya-u.ac.jp}

\urladdr{http://www.math.nagoya-u.ac.jp/\urltilda takahashi}

\thanks{Research partly supported by NSA grant H98230-11-0214 (L.W.C),
  and by JSPS Grant-in-Aid for Young Scientists (B) 22740008 and JSPS
  Postdoctoral Fellowships for Research Abroad (R.T)}

\date{16 November 2014}



% 13C05 Structure, classification theorems 13C60 Module categories
% 13C99 None of the above, but in this section 13E15 Rings and modules
% of finite generation or presentation; number of generators

\begin{abstract}
  A finitely generated module over a commutative noetherian ring of
  finite Krull dimension can be built from the prime ideals in the
  singular locus by iteration of three procedures: taking extensions,
  direct summands, and cosyzygies. In 2003 Schoutens gave a bound on
  the number of iterations required to build any module, and in this
  note we determine the exact number. This building process yields a
  stratification of the module category, which we study in detail for
  local rings that have an isolated singularity.
\end{abstract}

\ \vspace{-.25\baselineskip}

\maketitle

\thispagestyle{empty}

%%% INTRODUCTION

\section*{Introduction}

\noindent
Let $R$ be a commutative noetherian ring of finite Krull dimension. In
\cite{HSc03} Schoutens shows that starting from the set of singular
primes in $R$, one can build the entire category of finitely generated
$R$-modules by way of extensions, direct summands, and
cosyzygies. Schoutens's result gives a bound, in terms of the Krull
dimension of $R$, on the number of times these procedures must be
repeated to complete the building process. In this paper we give an
improved bound on this number and show that it is sharp. In the
process we give a condensed proof of the original~result.

From the building process one gets a stratification of the module
category into full subcategories that we call ``tiers''. Over a
regular ring the tiers simply sort the modules by projective
dimension, but over singular rings the picture remains opaque. We
describe the tiers explicitly for a local ring with an isolated
singularity.

\section{Tiers of modules}

\noindent
In this paper $R$ is a commutative noetherian ring, and $\mod{R}$
denotes the category of finitely generated $R$-modules. By a
subcategory of $\mod{R}$ we always mean a full subcategory closed
under isomorphisms. By $\RegR$ we denote the \emph{regular locus} of $R$;
that is, the set $\RegR = \setof{\p\in\Spec{R}}{\text{$R_\p$ is
    regular}}$. The \emph{singular locus} of $R$ is the complementary
set $\SingR = \Spec{R} \setminus \RegR$.

\begin{dfn}
  Let $\mathsf{S}$ be a  subcategory of $\mod{R}$.
  \begin{enumerate}[$\bullet$]
  \item Denote by $\close{\mathsf{S}}$ the smallest  subcategory
    of $\mod{R}$ that contains $\mathsf{S} \cup \set{0}$ and is closed
    under extensions and direct summands.
  \item Denote by $\cos{\mathsf{S}}$ the  subcategory whose
    objects are modules $X$ such that there exists an exact sequence
    $0 \to S \to P \to X \to 0$ where $S$ is in $\mathsf{S}$ and $P$
    is finitely generated and projective.
  \item Set $\tier{-1}{\mathsf{S}} = \close{\mathsf{S}}$,
    $\tier{0}{\mathsf{S}} = \close{\mathsf{S} \cup
      \cos{\close{\mathsf{S}}}}$ and for $n\in\NN$ set
    \begin{equation*}
      \tier{n}{\mathsf{S}} =
      \close{\tier{n-1}{\mathsf{S}} \cup
        \cos{(\tier{n-1}{\mathsf{S}})}}.
    \end{equation*}
  \end{enumerate}
\end{dfn}

\begin{ipg*}
  Let $\sing$ be the  subcategory of $\mod{R}$ with skeleton
  $\setof{R/\p}{\p\in\SingR}$; we consider the question of which, if
  any, of the subcategories in the chain
  \begin{equation*}
    \close{\sing} \: = \: \tier{-1}{\sing} \: \subseteq  \cdots
    \subseteq \: \tier{n}{\sing} \: \subseteq \:
    \tier{n+1}{\sing} \: \subseteq  \cdots
  \end{equation*}
  is the entire module category $\mod{R}$. In the rest of the paper, a
  subcategory of $\mod{R}$ described as a set $X$ is tacitly
  understood to be the subcategory with skeleton~$X$.

  In terms of of tiers, Schoutens's result
  \cite[Theorem~VI.8]{HSc03} can be stated as follows. If $R$ has
  finite Krull dimension $d$, then one has $\tier{d}{\sing} =
  \mod{R}$, and if $R$ is local and singular, then one has
  $\tier{d-1}{\sing}=\mod{R}$. For regular rings, Schoutens's bound is
  the best possible. Our theorem below sharpens the bound for singular
  rings: We replace $d$ (in the local case $d-1$) by
  $c=\codim{\mspace{2mu}(\SingR)}$, the \emph{codimension} of the
  singular locus, which is $-1$ if $\RegR$ is empty and otherwise
  given by
  \begin{equation*}
    c = \sup\setof{\hgt{\p}}{\p\in\RegR}.
  \end{equation*}
\end{ipg*}

\begin{thm} \label{main_thm} Let $R$ be a commutative noetherian ring
  and set $$\sing = \setof{R/\p}{\p\in\SingR}.$$ If\, $c =
  \codim{(\SingR)}$ is finite, then there is an equality
  $\tier{c}{\sing} = \mod{R}$.
\end{thm}

\enlargethispage*{\baselineskip}
\begin{prf*}
  As every $R$-module has a prime filtration and tiers are closed
  under extensions, it is sufficient to prove that every cyclic module
  $R/\p$, where $\p$ is a prime ideal in $R$, is in
  $\tier{c}{\sing}$. For a prime ideal $\p\in\RegR$, set
  \begin{equation*}
    n(\p) = \max\setof{\dim{(\q/\p)}}{\text{$\p \subseteq \q$ and $\q$
        is minimal in $\SingR$}}.
  \end{equation*}
  For $\p\in\SingR$, set $n(\p)=0$; we proceed by induction on
  $n(\p)$. By definition, $R/\p$ is in $\sing$ and, therefore, in
  $\tier{c}{\sing}$ if $n(\p)$ is $0$. Let $n\ge 1$ and assume that
  $R/\p$ is in $\tier{c}{\sing}$ for all $\p$ with $n(\p) < n$. Fix a
  prime ideal $\p$ with $n(\p)=n$ and set $h = \hgt{\p}$. Since $R_\p$
  is regular, one can choose elements $x_1,\dots, x_h$ in $\p$ such
  that the ideal $I = (x_1,\dots,x_h)$ has height $h$ and the equality
  \begin{equation}
    \label{eq:Ip}
    IR_\p = \p R_\p
  \end{equation}
  holds. As $\p/I$ is a minimal prime ideal in $R/I$, there exists an
  element $a \in R$ with $\p = (I\colon a)$, and it follows from
  \eqref{Ip} that $a$ is not in $\p$. It is now elementary to verify
  the equality $I = (I + (a))\cap \p$, which yields a Mayer--Vietoris
  exact sequence
  \begin{equation}
    \label{eq:mv}
    0 \lra R/I \lra R/\p \oplus R/(I + (a)) \lra R/(\p + (a)) \lra 0.
  \end{equation}
  The support of the module $R/(\p + (a))$ consists of prime ideals
  that strictly contain~$\p$. Thus, $R/(\p + (a))$ has a prime
  filtration with subquotients of the form $R/\q$, where each $\q$
  satisfies the inequality $n(\q) < n(\p)$. By the induction
  hypothesis, these subquotients $R/\q$ are in $\tier{c}{\sing}$ and
  hence so is $R/(\p + (a))$.

  By \eqref{mv} it now suffices to show that $R/I$ is in
  $\tier{c}{\sing}$.  To this end, consider the Koszul complex
  $K=K(x_1,\dots, x_h)$ on the generators of $I$. For $\q\in\RegR$,
  the non-units among the elements $x_1/1,\ldots,x_h/1$ in $R_\q$ form
  a regular sequence. It follows that the homology modules $\H[i]{K}$
  for $i > 0$ have support in $\SingR$, see \cite[Theorem~16.5]{mat},
  and therefore that they are in $\tier{-1}{\sing}$. Let
  $d_1,\ldots,d_h$ denote the differential maps on $K$. The modules
  $K_i$ in the Koszul complex are free, and the module $\Ker{d_h} =
  \H[h]{K}$ is in $\tier{-1}{\sing}$. It now follows from the exact
  sequences
  \begin{gather*}
    0 \lra \Im{d_{i+1}} \lra \Ker{d_i} \lra \H[i]{K} \lra 0\\
    0 \lra \Ker{d_i} \lra K_i \lra \Im{d_i} \lra 0
  \end{gather*}
  that $\Im{d_i}$ is in $\tier{h-i}{\sing}$ for $h \ge i \ge1$. In
  particular, the ideal $I=\Im{d_1}$ is in $\tier{h-1}{\sing}$. Thus
  the cosyzygy $R/I$ is in $\tier{h}{\sing}$ and clearly one has $h
  \le c$.
\end{prf*}

The proof above is quite close to Schoutens's original
argument.

\begin{rmk}
  One cannot leave out of any of the three procedures---adding
  cosyzygies, closing up under extensions, or closing up under
  summands---from the definition of tiers and still generate the
  entire module category. For the sake of the argument, let $R$ be an
  isolated curve singularity, i.e.\ a one-dimensional Cohen--Macaulay
  local ring $R$ with $\sing = \{k\}$, where $k$ is the residue field
  of $R$.
  \begin{enumerate}[$\bullet$]
  \item Without adding cosyzygies, one does not move beyond the
    category $\close{\sing}$, which contains only the $R$-modules of
    finite length and hence not $R$.
  \item The $R$-module $k$ is simple and cannot be embedded in a free
    $R$-module. Furthermore, $R$ is indecomposable as an
    $R$-module. It follows that by adding cosyzygies and closing up
    under summands one only gets $k$ and modules of projective
    dimension at most $1$. Thus, extensions are needed.
  \item Summands cannot be dispensed with either. The closure
    $\mathsf{E}$ of $\sing$ under extensions is the subcategory of
    modules of finite length. Since no such module can be embedded in
    a free $R$-module, $\cos{\mathsf{E}}$ contains exactly the
    finitely generated free modules. If the closure under extensions
    of $\mathsf{E} \cup \cos{\mathsf{E}}$ is the entire module
    category $\mod{R}$---or if $\mod{R}$ can be attained by
    alternately closing up under extensions and taking syzygies a
    finite number of times---then the Grothendieck group of $R$ is
    generated by $k$ and $R$. However for any even integer $n\ge 4$,
    the Grothendieck group of the $D_n$ singularity,
    $\mathbb{C}[\![x,y]\!]/(x^2y + y^{n-1})$, requires three
    generators; see \cite[Lemma~(13.2) and Proposition~(13.10)]{yos}.
  \end{enumerate}
\end{rmk}

\section{The codimension of $\SingR$ is the best possible bound}
\label{sec:result}

\noindent
We now show that the bound provided by Theorem \ref{main_thm} is
optimal; that is, $\tier{n}{\sing}$  for $n<c$ is a proper subcategory of
$\mod{R}$. First note that if $R$ is regular, then $\SingR $
and hence $\sing$ is empty. Thus $\tier{-1}{\sing}$ contains only the
zero module, and it follows from the definition that $\tier{n}{\sing}$
for $n\ge 0$ contains precisely the modules of projective dimension at
most $n$.  The next lemma shows that, to some extent, this simple
observation carries over to general rings.

\begin{lem}
  \label{lem:fl-pd}
  For a finitely generated $R$-module $M$ the following assertions
  hold.
  \begin{enumerate}[\rm(a)]
  \item $M$ is in $\tier{-1}{\sing}$ if and only if one has $M_\p=0$
    for every $\p\in\RegR$.
  \item If $M$ is in $\tier{n}{\sing}$, then
    $\pd[R_\p]{M_\p} \le n$ holds for every $\p\in\RegR$.
  \end{enumerate}
\end{lem}

\begin{prf*}
  As $\SingR$ is a specialization closed subset of $\Spec{R}$, one has
  $(R/\q)_\p=0$ for every $\q \in \SingR$ and every $\p\in\RegR$. It
  follows that $M_\p$ is $0$ for every $M\in \tier{-1}{\sing}$ and
  every $\p\in\RegR$. Conversely, if one has $M_\p=0$ for every
  $\p\in\RegR$, then $M$ has a prime filtration with subquotients
  $R/\q$ in $\sing$, so $M$ is in $\tier{-1}{\sing}$. This proves part
  (a).

  (b): Assume that $X$ is in $\cos{(\tier{-1}{\sing})}$, then there is
  exact sequence
  \begin{equation*}
    0 \lra S \lra P \lra X \lra 0,
  \end{equation*}
  where $P$ is a finitely generated projective module and $S$ is in
  $\tier{-1}{\sing}$. It follows that $X$ is free at every
  $\p\in\RegR$, and hence so are all modules in $\tier{0}{\sing}$.

  Let $n\ge 1$ and assume that the inequality $\pd[R_\p]{X_\p} \le
  n-1$ holds for all modules $X$ in $\tier{n-1}{\sing}$ and for every
  $\p\in\RegR$. It follows that every module in
  $\cos{(\tier{n-1}{\sing})}$ has projective dimension at most $n$ at
  every $\p\in\RegR$, and hence the desired inequality holds for all
  modules in $\tier{n}{\sing}$.
\end{prf*}

Up to $\tier{c}{\sing}$ each tier strictly contains the previous one.

\begin{prp}
  If\, $c = \codim{(\SingR)}$ is finite, then there are strict
  inclusions
  \begin{equation*}
    \tier{-1}{\sing} \:\subset\: \tier{0}{\sing} \:\subset \cdots \subset\:
    \tier{c-1}{\sing}  \:\subset\: \tier{c}{\sing}
  \end{equation*}
  of subcategories of $\mod{R}$.
\end{prp}

\begin{prf*}
  Let $\mathsf{S}$ be any subcategory of $\mod{R}$; if
  one has $\tier{n}{\mathsf{S}} = \tier{n+1}{\mathsf{S}}$ for some
  $n\ge -1$, then it follows from the definition that
  $\tier{n}{\mathsf{S}}$ equals $\tier{m}{\mathsf{S}}$ for all $m\ge
  n$.

  Thus, it is sufficient to show that $\tier{c-1}{\sing}$ is not the
  entire category $\mod{R}$. To this end choose a prime ideal $\p$ in
  $\RegR$ of height $c$. By the Auslander--Buchsbaum Equality one has
  $\pd[R_\p]{(R/\p)_\p} = c$, so it follows from Lemma \ref{lem:fl-pd}
  that $R/\p$ does not belong to $\tier{c-1}{\sing}$.
\end{prf*}

Our proof of Theorem~\ref{main_thm} only shows that every finitely
generated $R$-module is in $\tier{c}{\sing}$; it gives no
information on the least tier to which a given module $M$
belongs, but Lemma \ref{lem:fl-pd} provides a lower bound, namely
$\sup\setof{\pd[R_\p]{M_\p}}{\p\in\RegR}$.

Recall that a module $M \in \mod{R}$ is called is \emph{maximal
  Cohen--Macaulay} if the equality $\depth_R M = \dim R$ holds. Such a
module $M$ is free on the regular locus; indeed, the
Auslander--Buchsbaum Equality yields $\pd[R_\p]{M_\p}\le 0$ for all
$\p$ in $\RegR$. We show in the next section that over certain
Cohen--Macaualay local rings $R$ there are maximal Cohen--Macaulay
modules which are not in $\tier{0}{\sing}$.  Thus, the lower bound
provided by Lemma~\ref{lem:fl-pd} is not sharp, and we ask the
question:

\begin{bfhpg}[Question]\label{7}
  Let $R$ be a Cohen--Macaulay local ring and denote by $\cm(R)$ the
   subcategory of $\mod R$ consisting of all maximal
  Cohen--Macaulay modules.  What is the following number?
  \begin{equation*}
    \varepsilon(R)=\min\setof{n\ge -1}{\cm(R)\subseteq\tier n\s(R)}.
  \end{equation*}
\end{bfhpg}

If $R$ is a regular local ring, then $\varepsilon(R)$ is $0$ and we
show in the next section that it may be as big as $c=\codim{(\SingR)}$
for a singular ring. A broader question is, of course, given a module,
how can one determine the least tier it belongs to?

\section{Isolated singularities}

\noindent 
A local ring $R$ is Cohen--Macaulay if $R$ is a maximal
Cohen--Macaulay $R$-module, and $R$ is said to have an \emph{isolated
  singularity} if $R$ is singular but $R_\p$ is regular for every
non-maximal prime ideal in $R$. In this section we give a description
of the subcategories $\tier{n}{\sing}$ for a local ring $R$ with an
isolated singularity; one that is explicit enough to answer Question
\ref{7} for a Cohen--Macaualy local ring with an isolated singularity.

For a subcategory $\mathsf{S}$  of $\mod{R}$, every module in
$\close{\mathsf{S}}$ can be reached by alternate\-ly taking summands and
extensions; to discuss this we recall some notation from \cite{HLDRTk}.

\begin{dfn}
  Let $\S$ and $\T$ be subcategories of $\mod R$.
  \begin{enumerate}[(1)]
  \item Denote by $\add{\S}$ the additive closure of $\S$, that is,
    the smallest  subcategory of $\mod R$ containing $\S$ and
    closed under finite direct sums and direct summands.
  \item Denote by $\S\circ\T$ the subcategory of $\mod R$ consisting
    of the $R$-modules $M$ that fit into an exact sequence $ 0 \to S
    \to M \to T \to 0 $ with $S\in\S$ and $T\in\T$.
  \item Set $\S\bullet\T=\add{(\add{\S}\circ \add{\T})}$, and for
    integers $m\ge1$, set
    \begin{equation*}
      |\S|_m=
      \begin{cases}
        \add{\S} & \text{ for } m=1,\\
        |\S|_{m-1}\bullet\S & \text{ for } m\ge 2.
      \end{cases}
    \end{equation*}
  \end{enumerate}
\end{dfn}

\begin{rmk}
  Let $\S$ and $\T$ be subcategories of $\mod R$.  A module $M$ in
  $\mod R$ belongs to $\S\bullet\T$ if and only if there is an exact
  sequence $0 \to S \to E \to T \to 0$ with $S\in\add{\S}$ and $T\in\add{\T}$
  such that $M$ is a direct summand of $E$.  Moreover, one has
%  $(\S\bullet\T)\bullet\U=\S\bullet(\T\bullet\U)$ and
  $|\S|_m\bullet|\S|_{m'}=|\S|_{m+m'}$ for all $m,m' \ge 1$; see~\cite{HLDRTk}.
\end{rmk}

\begin{lem}
  \label{3}
  For every subcategory $\S$ of $\mod R$ one has
  $\close{\S}=\bigcup_{m\ge1}|\S|_m$ .
\end{lem}

\begin{prf*}
  Set $\T=\bigcup_{m\ge1}|\S|_m$.  Evidently one has
  $\S\subseteq\T\subseteq\close{\S}$, and $\T$ is by construction
  closed under direct summands.  Let
  \begin{equation*}
    0 \lra T \lra E \lra T' \lra 0
  \end{equation*}
  be an exact sequence in $\mod R$ with $T$ and $T'$ in $\T$.  There are
  integers $m,m' \ge 1$ with $T\in|\S|_{m}$ and $T' \in|\S|_{m'}$, and
  hence $E$ is in $|\S|_{m+m'}$.  Thus, $\T$ is also closed under
  extensions, and by the definition of $\close{S}$ it follows that one
  has $T=\close{S}$.
\end{prf*}

Let $R$ be a local ring with residue field $k$.  Denote by $\fl(R)$ the
 subcategory of $\mod{R}$ whose objects are all modules of finite
length. For $n\ge-1$ denote by $\PD_n(R)$ the  subcategory of
$\mod{R}$ whose objects are all modules of projective dimension at
most $n$. Note that one has $\fl(R) = \close{\set{k}}$ and
$\PD_{-1}(R)=\{0\}$.

\begin{thm}
  \label{5}
  Let $R$ be a local ring with residue field $k$.  For $-1\le
  n\le\depth R-1$ there are equalities of  subcategories of
  $\mod{R}$,
  \begin{equation*}
    \tier{n}{\set{k}} = \close{\fl(R)\cup\PD_n(R)} 
    = \close{\{k\}\cup\PD_n(R)},
  \end{equation*}
  and for $-1\le n\le\depth R-2$ the category $\tier{n}{\set{k}}$
  contains precisely the modules $M$ such that there is an exact
  sequence
  \begin{equation*}
    0 \lra L \lra M\oplus M' \lra P \lra 0
  \end{equation*}
  in $\mod R$ with $L\in\fl(R)$ and $P\in\PD_n(R)$.
\end{thm}

\begin{prf*}
  First we show that every module in $\close{\fl(R)\cup\PD_n(R)}$ for
  $-1 \le n \le \depth R - 2$ fits in an exact sequence $0 \to L \to
  M\oplus M' \to P \to 0$ with $L\in\fl(R)$ and $P\in\PD_n(R)$. The
  assertion is trivial for $n=-1$, so let $0\le n\le\depth R-2$.  Fix
  a module $M$ in $\close{\fl(R)\cup\PD_n(R)}$; by Lemma \ref{3} it
  belongs to $|\fl(R)\cup\PD_n(R)|_m$ for some $m\ge 1$. We now argue
  by induction on $m$ that $M$ fits in an exact sequence of the
  prescribed form.

  For $m=1$ one has $M\in\add{(\fl(R)\cup\PD_n(R))}$, whence there is
  an isomorphism $M\oplus M'\cong L\oplus P$ for modules $M'\in\mod
  R$, $L\in\fl(R)$, and $P\in\PD_n(R)$.

  For $m\ge2$ there is an exact sequence
  \begin{equation}
    \label{eq:4}
    0 \lra X \lra M \oplus M' \lra Y \lra 0
  \end{equation}
  in $\mod R$ with $X \in |\fl(R)\cup\PD_n(R)|_{m-1}$ and $Y \in
  \add{(\fl(R)\cup\PD_n(R))}$.  The base and hypothesis of
  induction yield an isomorphism $Y\oplus Y' \cong L\oplus P$ and an
  exact sequence $0 \to L' \to X\oplus X' \to P' \to 0$, with $L$ and
  $L'$ in $\fl(R)$ and with $P$ and $P'$ in $\PD_n(R)$.  Combined with
  \eqref{4} they yield an exact sequence
  \begin{equation*}
    0 \lra X\oplus X' \lra 
    X'\oplus M\oplus M' \oplus Y' \lra L\oplus P \lra 0.
  \end{equation*}
  Set $V=X'\oplus M \oplus M'\oplus Y'$. Consider the pushout diagram
  \begin{equation}
    \label{eq:po}
    \begin{gathered}
      \xymatrix@=1.73pc{ & 0 \ar[d] & 0 \ar[d] \\
        & L' \ar[d] \ar@{=}[r] & L' \ar[d] \\
        0 \ar[r] & X\oplus X' \ar[r] \ar[d] & V \ar[r] \ar[d] & L\oplus P \ar[r] \ar@{=}[d]& 0 \\
        0 \ar[r] & P' \ar[r] \ar[d] & W \ar[r] \ar[d] & L\oplus P \ar[r] & 0 \\
        & 0 & 0 \\ }
    \end{gathered}
  \end{equation}
  and the pullback diagram
  \begin{equation}
    \label{eq:pb}
    \begin{gathered}
      \xymatrix@=1.73pc{ && 0 \ar[d] & 0 \ar[d] \\
        0 \ar[r] & P' \ar[r] \ar@{=}[d] & P'' \ar[r] \ar[d] & P \ar[r] \ar[d] & 0 \\
        0 \ar[r] & P' \ar[r] & W \ar[r] \ar[d] & L\oplus P \ar[r] \ar[d] & 0 \\
        && L \ar@{=}[r] \ar[d] & L \ar[d] \\
        && 0 & 0 }
    \end{gathered}
  \end{equation}
  Note from the top row in \eqref{pb} that the module $P''$ is in
  $\PD_n(R)$. From the inequality $n \le \depth R - 2$ and the
  Auslander--Buchsbaum Equality one gets \mbox{$\depth_R P''\ge2$}. By
  the cohomological characterization of depth \cite[Theorem 16.6]{mat}
  this implies $\Ext_R^{1}(k,P'')=0$ and, therefore,
  $\Ext_R^1(L,P'')=0$. Thus, the middle column in \eqref{pb} is split
  exact, and the middle column in \eqref{po} becomes $0 \to L' \to V
  \to L\oplus P'' \to 0$.  Consider the pullback diagram
  \begin{equation}
    \label{eq:pb2}
    \begin{gathered}
      \xymatrix@=1.73pc{ && 0 \ar[d] & 0 \ar[d] \\
        0 \ar[r] & L' \ar[r] \ar@{=}[d] & L'' \ar[r] \ar[d] & L \ar[r] \ar[d] & 0 \\
        0 \ar[r] & L' \ar[r] & V \ar[r] \ar[d] & L\oplus P'' \ar[r] \ar[d] & 0 \\
        && P'' \ar@{=}[r] \ar[d] & P'' \ar[d] \\
        && 0 & 0 }
    \end{gathered}
  \end{equation}
  Note from the top row that $L''$ is in $\fl(R)$.  As $M$ is a direct
  summand of $V$, the middle column is a desired exact sequence.

  Clearly, one has $\close{\fl(R)\cup\PD_n(R)} =
  \close{\{k\}\cup\PD_n(R)}$; to finish the proof we show by
  induction that $\tier{n}{\set{k}} =
  \close{\fl(R)\cup\PD_n(R)}$ holds for $-1 \le n \le \depth R -1$. For
  $n=-1$, one has $\tier{n}{\set{k}}=\fl(R)=\close{\fl(R)\cup
    \PD_n(R)}$.  Let $n\ge0$ and assume that
  $\tier{n-1}{\set{k}}=\close{\fl(R)\cup\PD_{n-1}(R)}$ holds. By
  definition one then has
  \begin{equation*}
    \tier{n}{\set{k}} = \close{\close{\fl(R)\cup\PD_{n-1}(R)} 
      \cup \cosyz\close{\fl(R)\cup\PD_{n-1}(R)}},
  \end{equation*}
  whence it suffices to establish the equality
  \begin{equation*}
    \cosyz\close{\fl(R)\cup\PD_{n-1}(R)} = \PD_n(R).
  \end{equation*}
  The inclusion ``$\supseteq$'' is clear because a module in
  $\PD_n(R)$ is a cosyzygy of its first syzygy, which is in
  $\PD_{n-1}(R)$. For the opposite inclusion, let $M$ be a module in
  $\cosyz\close{\fl(R)\cup\PD_{n-1}(R)}$.  There is an exact
  sequence
  \begin{equation*}
    0 \lra N \lra F \lra M \lra 0,
  \end{equation*}
  where $F$ is free and $N$ is in $\close{\fl(R)\cup\PD_{n-1}(R)}$.
  From the inequalities $-1 < n \le \depth R - 1$ follows that $R$ and
  hence $F$ has positive depth, whence also $N$ has positive depth.
  Moreover, one has $-1\le n-1\le\depth R-2$, so it follows from the
  first part of the proof that there is an exact sequence in
  $\mod{R}$,
  \begin{equation*}
    0 \lra L \xrightarrow{\binom{\alpha}{\beta}} N\oplus N' \lra P \lra 0,
  \end{equation*}
  with $L\in\fl(R)$ and $P\in\PD_{n-1}(R)$.  Since $L$ has finite
  length and $N$ has positive depth, the map $\alpha$ is zero. Thus,
  there is an isomorphism $P \cong N \oplus C$, where $C$ is the
  cokernel of $\beta$.  Thus $N$ belongs to $\PD_{n-1}(R)$, and hence
  $M$ is in $\PD_n(R)$.
\end{prf*}

For a local ring $R$ with an isolated singularity one has $\sing =
\set{k}$, so Theorems~\ref{main_thm} and \ref{5} combine to yield:

\begin{cor}    
  \protect\pushQED{\qed} Let $R$ be a $d$-dimensional local ring with
  an isolated singularity.  For every $n\ge -1$ one has
  \begin{equation*}
    \tier n \sing = \close{\set{k}\cup\PD_{n}(R)};
  \end{equation*}
  in particular, one has
  \begin{equation*}
    \mod{R} = \close{\set{k}\cup\PD_{d-1}(R)}.\qedhere
  \end{equation*}
\end{cor}

To answer Question \ref{7} for a Cohen--Macaualy local ring with an
isolated singularity, we record another consequence of Theorem
\ref{5}.  Denote by $\dep(R)$ the  subcategory of $\mod{R}$ whose
objects are all modules of positive depth; it includes the zero module
as it has infinite depth by convention.

\begin{prp}
  \label{8}
  Let $R$ be a local ring. For $-1 \le
  n\le\depth R-2$ one has
  \begin{equation*}
    \tier{n}{\set{k}}\cap\dep(R)=\PD_n(R).
  \end{equation*}
\end{prp}

\begin{prf*}
  By Theorem \ref{5} one has
  $\tier{n}{\set{k}}=\close{\fl(R)\cup\PD_n(R)}$ and it follows from
  the inequality $n \le \depth R -2$ and the Auslander--Buchsbaum
  Equality that every module in $\PD_n(R)$ has positive depth.  This
  proves the inclusion ``$\supseteq$''.  To show the opposite
  inclusion, fix a module $M$ in
  $\close{\fl(R)\cup\PD_n(R)}\cap\dep(R)$.  It follows from Theorem
  \ref{5} that there is an exact sequence in $\mod R$
  \begin{equation*}
    0 \lra L \xrightarrow{\binom{\alpha}{\beta}} M\oplus M' \lra P \lra 0
  \end{equation*}
  with $L\in\fl(R)$ and $P\in\PD_n(R)$.  Since $M$ has positive depth,
  the map $\alpha$ is zero and it follows that $M$ is a direct summand
  of $P$, whence $M$ is in $\PD_n(R)$.
\end{prf*}

\begin{cor}
  \label{cor:notsharp}
  Let $R$ be a $d$-dimensional Cohen--Macaulay local ring.
  \begin{enumerate}[\rm(a)]
  \item If $d\ge1$, then every maximal Cohen--Macaulay module in
    $\tier{d-2}{\set{k}}$ is free.
  \item If $R$ has an isolated singularity, then one has
    \begin{equation*}
      \varepsilon(R)=d-1.
    \end{equation*}
  \end{enumerate}
\end{cor}

\begin{prf*}
  (a) By Proposition \ref{8} one has
  \begin{equation*}
    \cm(R)\cap\tier{d-2}{\set{k}} = \cm(R)\cap(\tier{d-2}{\set{k}}
    \cap\dep(R))=\cm(R)\cap\PD_{d-2}(R).
  \end{equation*}
  In $\cm(R)\cap\PD_{d-2}(R)$ is only $0$ if $d=1$ and precisely the
  free $R$-modules if $d\ge2$.

  (b) The equality is trivial for $d=0$ and it follows from (a) for
  $d\ge1$.
\end{prf*}

The corollary shows that the lower bound that Lemma~\ref{lem:fl-pd}
gives for which tier a module $M$ can belong to,
$\sup\setof{\pd[R_\p]{M_\p}}{\p\in\RegR}$, is far from being sharp.

\bibliographystyle{amsplain} %\bibliography{../+references}
\providecommand{\MR}[1]{\mbox{\href{http://www.ams.org/mathscinet-getitem?mr=#1}{#1}}}
  \renewcommand{\MR}[1]{\mbox{\href{http://www.ams.org/mathscinet-getitem?mr=#1}{#1}}}
\providecommand{\bysame}{\leavevmode\hbox to3em{\hrulefill}\thinspace}
\providecommand{\MR}{\relax\ifhmode\unskip\space\fi MR }
\providecommand{\MRhref}[2]{%
  \href{http://www.ams.org/mathscinet-getitem?mr=#1}{#2}
}
\providecommand{\href}[2]{#2}

\end{document}